\documentclass[11pt]{article}
\usepackage{graphicx}
\usepackage{amssymb,latexsym}
\usepackage{epstopdf}
\newtheorem{lemma}{Lemma}[section]
\newtheorem{theorem}[lemma]{Theorem}

\newtheorem{claim}[lemma]{Claim}

\newtheorem{remark}[lemma]{Remark}

\newcommand{\drg}{distance-regular graph}
\newcommand{\srg}{strongly regular graph}
\newcommand{\g}{\Gamma}
\newcommand{\pf}{\noindent{\em Proof: }}
\newcommand{\epf}{\hfill\hbox{\rule{3pt}{6pt}}\\}
\newcommand{\forme}[1]{}
\newcommand{\HH}{{\texttt h}}

\begin{document}

\title{Geometric distance-regular graphs without $4$-claws}

\author{
{\bf S.~Bang }
 \\
Department of Mathematical Sciences, BK21-Mathematical Sciences Division\\
Seoul National University, Seoul 151-742 South Korea\\
e-mail: sjbang3@snu.ac.kr}

\date{\today}

\maketitle
\begin{abstract}
A non-complete \drg~$\Gamma$ is called geometric if there exists a
set $\mathcal{C}$ of Delsarte cliques such that each edge of
$\Gamma$ lies in a unique clique in $\mathcal{C}$. In this paper, we
determine the non-complete distance-regular graphs satisfying $\max
\{3, \frac{8}{3}(a_1+1)\}<k<4a_1+10-6c_2$. To prove this result, we
first show by considering non-existence of $4$-claws that any
non-complete distance-regular graph satisfying $\max \{3,
\frac{8}{3}(a_1+1)\}<k<4a_1+10-6c_2$ is a geometric \drg~with
smallest eigenvalue $-3$. Moreover, we classify the geometric \drg s
with smallest eigenvalue $-3$. As an application,  $7$ feasible
intersection arrays in the list of \cite[Chapter 14]{bcn} are ruled
out.
\end{abstract}

\section{Introduction}\label{intro}
Let $\Gamma$ be a \drg~with valency $k$ and let $\theta_{\min}=\theta_{\min}(\g)$ be its smallest eigenvalue. Any clique $C$ in $\Gamma$ satisfies
\begin{equation}\label{hoffman-bd}
|C|\leq 1-\frac{k}{\theta_{\min}}
\end{equation}
(see \cite[Proposition 4.4.6 (i)]{bcn}). This bound
(\ref{hoffman-bd}) is due to Delsarte, and a clique $C$ in $\Gamma$
is called a {\em Delsarte clique} if $C$ contains exactly
$1-\frac{k}{\theta_{\min}}$ vertices. Godsil \cite{godsil-93-paper}
introduced the following notion of a geometric \drg. A non-complete
\drg~$\Gamma$ is called {\em geometric} if there exists a set
$\mathcal{C}$ of Delsarte cliques such that each edge of
$\Gamma$ lies in a unique Delsarte clique in $\mathcal{C}$. In this case, we say that $\g$ is geometric with respect to $\mathcal{C}$.\\

There are many examples of geometric distance-regular graphs such as
bipartite distance-regular graphs, the Hamming graphs, the Johnson
graphs, the Grassmann graphs and regular near $2D$-gons.\\
In particular, the local structure of geometric distance-regular
graphs play an important role in the study of spectral
characterization of some distance-regular graphs. In \cite{H1}, we
show that for given integer $D\geq 2$, any graph cospectral with the
Hamming graph $H(D,q)$ is locally the disjoint union of $D$ copies
of the complete graph of size $q-1$, for $q$ large enough. By using
this result and \cite{H0}, we show in \cite{H1} that the Hamming
graph
$H(3,q)$ with $q\geq 36$ is uniquely determined by its spectrum. \\

Neumaier \cite{neumaier-m} showed that except for a finite number of
graphs, any geometric strongly regular graph with a given smallest
eigenvalue $-m$, $m>1$ integral,
is either a Latin square graph or a Steiner graph (see \cite{neumaier-m} and Remark \ref{geo-rmk} for the definitions). \\
An {\em $n$-claw} is an induced subgraph on $n+1$ vertices which consists of one vertex of valency
$n$ and $n$ vertices of valency $1$. Each \drg~without $2$-claws is a complete graph.
Note that for any geometric \drg~$\g$ with respect to $\mathcal{C}$ a set of Delsarte cliques, the number of Delsarte cliques in $\mathcal{C}$ containing a fixed vertex is $-\theta_{\min}(\g)$. Hence any geometric \drg~with smallest eigenvalue $-2$ contains no $3$-claws. Blokhuis and Brouwer \cite{3-claw} determined the \drg s without $3$-claws. \\
Yamazaki \cite{yamazaki} considered \drg s which are locally a disjoint union of three cliques of size $a_1+1$, and these graphs for $a_1\geq 1$ are geometric \drg s with smallest eigenvalue $-3$.\\

In Theorem \ref{gdrg}, we determine the geometric \drg s with
smallest eigenvalue $-3$. We now state our main result of this
paper.

\begin{theorem}\label{main-cor}
Let $\g$ be a non-complete distance-regular graph. If $\g$ satisfies \[ \max \{3, \frac{8}{3}(a_1+1)\}<k<4a_1+10-6c_2\] then $\g$ is one of the following.
\begin{enumerate}
\item[(i)] A Steiner graph $S_3(\alpha -3)$, i.e., a geometric strongly regular graph with parameters $\left(\frac{(2\alpha-3)(\alpha-2)}{3},3\alpha-9,\alpha,9  \right)$, where $\alpha\geq 36$ and $\alpha \equiv 0,2~~(\mbox{mod}~3)$.
\item[(ii)] A Latin square graph $LS_3(\alpha)$, i.e., a geometric strongly regular graph with parameters $(\alpha ^2,3(\alpha-1),\alpha,6)$, where $\alpha \geq 24$.
\item[(iii)] The generalized hexagon of order $(8,2)$ with $\iota(\Gamma)=\{24,16,16;1,1,3\}$.
\item[(iv)] One of the two generalized hexagons of order $(2,2)$ with $\iota(\g)=\{6,4,4;1,1,3\}$.
\item[(v)] A generalized octagon of order $(4,2)$ with $\iota(\g)=\{12,8,8,8;1,1,1,3\}$.
\item[(vi)] The Johnson graph $J(\alpha,3)$, where $\alpha \geq 20$.
\item[(vii)] $D=3$ and $\iota(\g)=\{3\alpha +3,2\alpha +2, \alpha +2-\beta ;1,2,3\beta\}$, where $\alpha \geq 6$ and $\alpha\geq \beta \geq 1$.
\item[(viii)] The halved Foster graph with $\iota(\g)=\{6,4,2,1;1,1,4,6\}$.
\item[(ix)] $D=\mbox{\em \HH}+2\geq 4$ and
$$(c_i,a_i,b_i)=\left\{ \begin{array}{ll}
(1,\alpha,2\alpha+2)& \mbox{ for }1\leq i\leq \mbox{\em \HH}\\
(2,2\alpha+\beta-1,\alpha-\beta+2) & \mbox{ for }i=\mbox{\em \HH} +1\\
(3\beta,3\alpha-3\beta+3,0) & \mbox{ for }i=\mbox{\em \HH} +2
\end{array}
\right., \mbox{~where~} \alpha \geq  \beta \geq 2.$$
\item[(x)] $D=\mbox{\em \HH}+2\geq 3$ and $$(c_i,a_i,b_i)=\left\{ \begin{array}{ll}
(1,\alpha,2\alpha+2)& \mbox{ for }1\leq i\leq \mbox{\em \HH} \\
(1,\alpha+2\beta-2,2\alpha-2\beta+4) & \mbox{ for }i=\mbox{\em \HH} +1\\
(3\beta,3\alpha-3\beta+3,0) & \mbox{ for }i=\mbox{\em \HH} +2
\end{array}
\right., \mbox{~where~} \alpha\geq \beta \geq 2.$$
\item[(xi)] A distance-$2$ graph of a distance-biregular graph with vertices of valency $3$ and
$$(c_i,a_i,b_i)=\left\{ \begin{array}{ll}
(1,\alpha,2\alpha+2)& \mbox{ for }1\leq i\leq \mbox{\em \HH}\\
(1,\alpha+2,2\alpha)& \mbox{ for }i=\mbox{\em \HH} +1\\
(4,2\alpha-1,\alpha)& \mbox{ for }\mbox{\em \HH} +2\leq i\leq D-2\\
(4,2\alpha+\beta-3,\alpha-\beta+2) & \mbox{ for }i=D-1\\
(3\beta,3\alpha-3\beta+3,0) & \mbox{ for }i=D
\end{array}
\right., \mbox{~where~} \alpha\geq \beta \mbox{~and~} \beta\in \{2,3\}.$$
\end{enumerate}
\end{theorem}

Examples of non-complete \drg s with valency $k>\max \{3, \frac{8}{3}(a_1+1)\} $ include Johnson graphs $J(n,e)$
$\Big( (n\geq 20 \mbox{~and~}e=3)$, $(n\geq 11 \mbox{~and~}e=4)$ or $(n\geq 2e \mbox{~and~}e\geq 5) \Big)$, Hamming graphs $H(d,q)$ \Big($(d=3
\mbox{~and~}q\geq 3)$~or~$(d\geq 4\mbox{~and~}q\geq 2 )$\Big) and
Grassmann graphs $\Bigl[{{V}\atop e}\Bigr]$ $\Big( (e=2\mbox{~and~}q\geq 4 )$~or~$(e\geq 3 \mbox{~and~}q\geq 2) \Big)$, where $n\geq 2e$ and $V$ is an $n$-dimensional vector space over $\mathbb{F}_q$ the finite field of $q$($\geq 2$) elements (see \cite[Chapter 9]{bcn} for more information on these examples). Except $J(n,3)~(n\geq 20)$ and $H(3,q)~(q\geq 3)$, all the above
examples contain $4$-claws. Whereas, $J(n,3)~(n\geq 20)$ and
$H(3,q)~(q\geq 3)$ are geometric \drg s with smallest eigenvalue $-3$.\\

In Section \ref{no-4-claws}, we prove Theorem \ref{main-thm} which gives a sufficient condition, $\max \{3, \frac{8}{3}(a_1+1)\}<k<4a_1+10-6c_2$, for geometric \drg s with smallest eigenvalue $-3$. We first show in Theorem \ref{geo} that for any \drg~satisfying $k>\max \{3, \frac{8}{3}(a_1+1)\}$, the statement that $\g$ has no $4$-claws is equivalent to the statement that $\g$ is geometric with smallest eigenvalue $-3$. By using Theorem \ref{geo}, we will prove Theorem \ref{main-thm}. As an application of Theorem \ref{geo}, we can show non-existence of a family of \drg s with feasible intersection arrays. For example, in the list of \cite[Chapter 14]{bcn}, the $7$ feasible intersection arrays in Theorem \ref{non-exist-7} are ruled out.\\

In Section \ref{geosection}, we determine the geometric \drg s with smallest eigenvalue $-3$ in Theorem \ref{gdrg}. By using Theorem \ref{main-thm} and Theorem \ref{gdrg}, we will prove Theorem \ref{main-cor}.

\section{Preliminaries}\label{pre}
All graphs considered in this paper are finite, undirected and
simple (for unexplained terminology and more details, see
\cite{bcn}).\\
For a connected graph $\Gamma$, distance $d_{\Gamma}(x,y)$
between any two vertices $x,y$ in the
vertex set $V(\Gamma)$ of $\Gamma$ is the length of a shortest
path between $x$ and $y$ in $\Gamma$, and denote by $D(\Gamma)$ the diameter of $\Gamma$ (i.e., the
maximum distance between any two vertices of $\Gamma$). For any vertex $x\in V(\Gamma)$, let $\Gamma_i(x)$ be the set of vertices in
$\Gamma$ at distance precisely $i$ from $x$, where $i$ is a
non-negative integer not exceeding $D(\Gamma)$. In addition, define
$\Gamma_{-1}(x)=\Gamma_{D(\Gamma)+1}(x):=\emptyset$ and $\Gamma_0(x):=\{x\}$. For any
distinct vertices $x_1,x_2,\ldots, x_j\in V(\Gamma)$, define
$$\Gamma_1(x_1,\ldots,x_j):=\Gamma_1(x_1)\cap \Gamma_1(x_2) \cap
\cdots \cap \Gamma_1(x_j).$$
A {\em clique} is a set of pairwise adjacent vertices. A graph $\g$ is called {\em locally} $G$ if any local graph of $\g$ (i.e., the local graph of a vertex $x$ is the induced subgraph on $\Gamma_1(x)$) is isomorphic to $G$, where $G$ is a graph. The {\em adjacency matrix} $A(\Gamma)$ of a graph $\Gamma$ is the
$|V(\Gamma)|\times |V(\Gamma)|$-matrix with rows and columns are
indexed by $V(\Gamma)$, and the $(x,y)$-entry of $A(\g)$ equals $1$
whenever $d_{\g}(x,y)=1$ and $0$ otherwise. The eigenvalues of $\Gamma$ are the eigenvalues of $A(\g)$.\\
A connected graph $\Gamma$ is called a {\em\drg}~if there exist integers $b_i(\g)$, $c_i(\g)$, $i=0,1,\ldots,D(\Gamma)$, such that for any
two vertices $x,y$ at distance $i=d_{\Gamma}(x,y)$, there are
precisely $c_i(\g)$ neighbors of $y$ in $\Gamma_{i-1}(x)$ and $b_i(\g)$
neighbors of $y$ in $\Gamma_{i+1}(x)$. In particular, $\Gamma$ is
regular with valency $k(\g):=b_0(\g)$. The numbers $c_i(\g),b_i(\g)$ and
$a_i(\g):=k(\g)-b_i(\g)-c_i(\g)~(0\leq i\leq D(\g))$ (i.e., the number of neighbors of
$y$ in $\Gamma_i(x)$ for $d_{\g}(x,y)=i$) are called the {\em
intersection numbers} of $\Gamma$. Note that $b_{D(\g)}(\g)=c_0(\g)=a_0(\g):=0$ and
$c_1(\g)=1$. In addition, we define $k_i(\g):=|\Gamma_i(x)|$ for any vertex
$x$ and $i=0,1,\ldots,D(\g)$. The array
$\iota(\g)=\{b_0(\g),b_1(\g),\ldots,b_{D(\g)-1}(\g);c_1(\g),c_2(\g),\ldots,c_{D(\g)}(\g)\}$ is called the {\em
intersection array} of $\Gamma$. In addition, we define the number
\begin{equation}\label{head}
\HH(\Gamma) := |\{j \,\mid \, (c_j,a_j,b_j)
=(c_1,a_1,b_1),\,1\le j \le D(\g)-1\}|
\end{equation}
which is called the {\em head} of $\Gamma$. \\
A regular graph $\Gamma$ on $v$ vertices with valency $k(\g)$ is called a {\em
strongly regular graph} with parameters $(v,k(\g),\lambda(\g),\mu(\g))$ if there
are two constants $\lambda(\g)\geq 0$ and $\mu(\g)>0$ such that for any two distinct
vertices $x$ and $y$, $|\Gamma_1(x,y)|$ equals $\lambda(\g)$ if $d_{\g}(x,y)=1$ and $\mu(\g)$ otherwise.\\
When there are no confusion, we omit $\sim_{\Gamma}$ and $\sim(\g)$ in each
notation for $\Gamma$, such as $d_{\Gamma}(~,~)$, $D(\g)$, $A(\g)$, $\HH(\g)$, $k(\Gamma), c_i(\Gamma), b_i(\Gamma), a_i(\Gamma)$, $k_i(\g)$, $\lambda(\g)$ and $\mu(\g)$.\\
Suppose that $\Gamma$ is a \drg~with valency $k\geq 2$ and diameter
$D\geq 2$. It is well-known that $\Gamma$ has
exactly $D+1$ distinct eigenvalues which are the eigenvalues of the
following tridiagonal matrix
\begin{equation}\label{mtx-L}
 L_1(\Gamma):= \left\lgroup
 \begin{tabular}{llllllll}
 $0$ & $b_0$\\
 $c_1$ & $a_1$ & $b_1$\\
 & $c_2$ & $a_2$ & $b_2$\\
 & & . & . & .\\
 & & & $c_i$ & $a_i$ & $b_i$\\
 & & & & . & . & .\\
 & & & & & $c_{D-1}$ & $a_{D-1}$ & $b_{D-1}$\\
& & & & &\makebox{\hspace{.324cm}} &
 $c_{D}$ & $a_{D}$
 \end{tabular}
 \right\rgroup
\end{equation}
(cf. \cite[p.128]{bcn}). In particular, we denote by $\theta_{\min}=\theta_{\min}(\g)$ the smallest eigenvalue of $\g$.

\section{Distance-regular graphs without $4$-claws}\label{no-4-claws}
In this section, we prove the following theorem which gives a sufficient condition for geometric \drg s with smallest eigenvalue $-3$.

\begin{theorem} \label{main-thm}
Let $\g$ be a non-complete distance-regular graph. If $\g$ satisfies
\begin{equation}\label{k-condi}
 \max \{3, \frac{8}{3}(a_1+1)\}<k<4a_1+10-6c_2
\end{equation}
 then $\g$ is a geometric \drg~with smallest eigenvalue $-3$.
\end{theorem}

We first show in Theorem \ref{geo} that for any \drg~satisfying $k>\max \{3, \frac{8}{3}(a_1+1)\}$, the statement that $\g$ has no $4$-claws is equivalent to the statement that $\g$ is geometric with smallest eigenvalue $-3$. By using Theorem \ref{geo}, we will prove Theorem \ref{main-thm}. As an application, by considering a restriction on $c_2$ in Lemma \ref{mu-bd}, we can rule out a family of feasible intersection arrays. In particular, we prove that there are no \drg s with the intersection arrays in Theorem \ref{non-exist-7}.

\begin{theorem}\label{geo}
Let $\Gamma$ be a \drg ~satisfying $k >
\max \{3, \frac{8}{3}(a_1+1)\}$. Then the following are equivalent.\\
(i) $\Gamma$ has no $4$-claws.\\
(ii) $\Gamma$ is a
geometric \drg~with smallest eigenvalue $-3$.
\end{theorem}

\pf Let $\Gamma$ be a \drg~satisfying $k>\max \{3, \frac{8}{3}(a_1+1)\}$. Let $\theta_{\min}=\theta_{\min}(\Gamma)$.\\
\noindent (ii)$\Rightarrow $(i): Suppose that $\g$ is geometric with respect to $\mathcal{C}$ a set of Delsarte cliques and $\theta_{\min}=-3$. Since the number of Delsarte cliques in $\mathcal{C}$ containing a given vertex is $-\theta_{\min}$, the statement (i) follows immediately.\\
\noindent (i)$\Rightarrow $(ii): Suppose that $\Gamma$ has no $4$-claws. Define a {\em line} to be a
maximal clique $C$ in $\Gamma$ such that $C$ has at least
$k-2(a_1+1)+1$ vertices. Note here that $a_1\geq 1$ follows,
otherwise $\Gamma$ has a $4$-claw from
$k>\max \{3, \frac{8}{3}(a_1+1)\}$. Hence, $|C|\geq 3$ for any line $C$ in $\Gamma$. If there exists a line $C$ satisfying $|C|=3$, then $a_1=1$ and $k=6$ both hold by $3\geq k-2(a_1+1)+1$ and $k > \frac{8}{3}(a_1+1)$. By
\cite[Theorem 1.1]{a=1k=6}, the graph $\g$ is one
of the following.\\
(a) The generalized quadrangle of order $(2,2)$.\\
(b) One of the two generalized hexagons of order $(2,2)$.\\
(c) The Hamming graph $H(3,3)$.\\
(d) The halved Foster graph.\\
All the graphs in (a)-(d) are geometric with smallest eigenvalue $-3$.\\
In the rest of the proof, we assume that each line contains more
than $3$ vertices. First, we prove the following claim.
\begin{claim}\label{geo-claim1}
Every edge of $\Gamma$ lies in a unique line.
\end{claim}
\noindent{\em Proof of Claim \ref{geo-claim1}: } Let $(x,y_1)$ be an
arbitrary edge in $\Gamma$. As $k\geq 2(a_1+1)+1$, there exists a
$3$-claw containing $x$ and $y_1$, say $\{x,y_1,y_2,y_3\}$ induces a $3$-claw, where
$y_i\in \Gamma_1(x)$ ($i=1,2,3$). Put $Y_i:=\{y_i\}\cup
\Gamma_1(x,y_i)$ ($i=1,2,3$). If there exists a vertex $z$ in $\Gamma_1(x)\setminus \cup_{i=1}^{3}Y_i$, then $\{x,z,y_1,y_2,y_3\}$
induces a $4$-claw which is impossible, and therefore $\Gamma_1(x)=\cup_{i=1}^{3}Y_i$
follows. If there exist non-adjacent two vertices $v,w$ in
$Y_1\setminus (Y_2\cup Y_3)$, then the set $\{x,y_2,y_3,v,w\}$
induces a $4$-claw which is a contradiction. Hence $\{x\}\cup
\left(Y_1\setminus (Y_2\cup Y_3)\right) $ induces a clique
containing the edge $(x,y_1)$, and it satisfies
\[\left| \{x\}\cup \left(Y_1\setminus \left(Y_2\cup Y_3\right)\right) \right|
=|\{x\}|+|\Gamma_1(x)|-|Y_2\cup Y_3|\geq 1+k-2(a_1+1). \] Thus
every edge lies in a line.\\
Assume that there exist two lines $C_z$ and $C_w$ containing the edge
$(x,y_1)$, where $z\in C_z$ and $w\in C_w$ are two non-adjacent vertices. Then $a_1=|\Gamma_1(x,y_1)|\geq 2(k-2(a_1+1)-1)-(|C_z\cap C_w|-2)$ implies
\begin{equation}\label{two lines t}
|C_z\cap C_w|\geq 2k-5a_1-4.
\end{equation}
In addition, by (\ref{two lines t}),
\begin{eqnarray}\label{k-bd(geo)}
\left| \Gamma_1(x)\setminus \left(\Gamma_1(x,z)\cup
\Gamma_1(x,w)\cup\{z,w\}\right)\right|&\geq &
k-(\left|\Gamma_1(x,z)\right|
+\left|\Gamma_1(x,w) \right|+\left| \{z,w\}\right|-(|C_z\cap C_w|-1))\nonumber\\
&\geq & k-(2(a_1+1)-(2k-5a_1-5)) \nonumber \\&=& 3k-7a_1-7.
\end{eqnarray}
Since $\Gamma$ has no $4$-claws, $(\{x\}\cup \Gamma_1(x))\setminus
\left(\Gamma_1(x,z)\cup \Gamma_1(x,w)\cup\{z,w\}\right)$ induces a
clique of size at least $3k-7a_1-6$ by (\ref{k-bd(geo)}). Since any
clique in $\g$ has size at most $a_1+2$, we have $k\leq \frac{8}{3}(a_1+1)$
which is impossible. Hence, the edge $(x,y_1)$ lies in a unique line. Now, Claim \ref{geo-claim1} is proved. \epf

For each vertex $x\in V(\Gamma)$, we define $M_x$ to be the number
of lines containing $x$. Then for any vertex $x$, we have $M_x\geq 3$ as $k>\frac{8}{3}(a_1+1)>2(a_1+1)$, and hence
\begin{equation}\label{mx=3}
M_x=3 \mbox{ for each vertex } x\in V(\Gamma)
\end{equation}
as $k\geq M_x(k-2(a_1+1))$ holds by Claim \ref{geo-claim1}. Let $B$ be the vertex-line incidence matrix (i.e., the
$(0,1)$-matrix with rows and columns are indexed by the vertex set
and the set of lines of $\Gamma$ respectively, where $(x,C)$-entry
of $B$ is $1$ if the vertex $x$ is contained in the line $C$ and $0$
otherwise). By Claim \ref{geo-claim1} and (\ref{mx=3}),
$BB^T=A+3I$ holds, where $B^T$ is the transpose of $B$, $A=A(\g)$ and
$I$ is the $|V(\g)|\times |V(\g)|$ identity matrix. Since each line contains more than $3$
vertices, it follows by double-counting the number of ones in $B$ that the number of lines is strictly less than the number of vertices in $\Gamma$. Hence, the matrix $BB^T$ is singular so that $0$ is an eigenvalue of $BB^T$ and thus $-3$ is an eigenvalue of $A$. As $BB^T$ is positive semidefinite, we find $\theta_{\min}=-3$. Hence it follows by (\ref{hoffman-bd}), Claim \ref{geo-claim1}, (\ref{mx=3}) and $\theta_{\min}=-3$ that every line has exactly
$1+\frac{k}{3}$ vertices. This proves that $\Gamma$ is geometric with $\theta_{\min}=-3$. \epf

In \cite[Lemma 2]{shilla}, Koolen and Park have shown the following lemma.

\begin{lemma} \label{mu-bd}
Let $\Gamma$ be a \drg~with a $4$-claw. Then $\g$ satisfies
\[c_2\geq \frac{4a_1+10-k}{6}.\]
\end{lemma}
\pf Suppose that $\{x,y_i\mid 1\leq i\leq 4\}$ induces a $4$-claw in $\g$,
where $y_i\in \Gamma_1(x)$ ($i=1,2,3,4$). It follows by the principle of
inclusion and exclusion that
\begin{eqnarray*}
k &\geq & \left|\{y_i\mid 1\leq i\leq 4\} \right|+\left|
\cup_{i=1}^{4} \Gamma_1(x,y_i)\right|\nonumber \\ &\geq &
\left|\{y_i\mid 1\leq i\leq 4\} \right|+\sum_{i=1}^{4}
\left|\Gamma_1(x,y_i)\right|-\sum_{1\leq i< j\leq
4}\left|\Gamma_1(x,y_i,y_j)\right|\nonumber \\ &\geq & 4+4a_1-
{4\choose 2}  (c_2-1),
\end{eqnarray*}
from which Lemma \ref{mu-bd} follows. \epf

We now prove our main result of Section \ref{no-4-claws}, Theorem \ref{main-thm}.\\

\noindent{\em Proof of Theorem \ref{main-thm}:} Suppose that $\g$ is a non-complete distance-regular graph satisfying (\ref{k-condi}). Then there are no $4$-claws in $\g$ by Lemma \ref{mu-bd}, so that $\g$ is geometric with $\theta_{\min}(\g)=-3$ by Theorem \ref{geo}. This completes the proof. \epf

\begin{theorem}\label{non-exist-7}
There are no \drg s with the following intersection arrays\\
(i) $\{55,36,11;1,4,45\}$,\\
(ii) $\{56,36,9;1,3,48\}$,\\
(iii) $\{65,44,11;1,4,55\}$,\\
(iv) $\{81,56,24,1;1,3,56,81\}$,\\
(v) $\{117,80,32,1;1,4,80,117\}$,\\
(vi) $\{117,80,30,1;1,6,80,117\}$,\\
(vii) $\{189,128,45,1;1,9,128,189\}$.
\end{theorem}

\pf Assume that $\Gamma$ is a \drg~such that its intersection array is one of the $7$ intersection arrays (i)-(vii). Since $\Gamma$ satisfies $k>
\frac{8}{3}(a_1+1)$, $a_1\neq 0$ and $\theta_{\min}(\Gamma)\neq -3$, $\Gamma$ has a $4$-claw by Theorem \ref{geo}. It follows by Lemma \ref{mu-bd} that $c_2\geq \frac{4a_1+10-k}{6}$ which is impossible. This shows Theorem \ref{non-exist-7}.\epf

\begin{remark}
\begin{enumerate}
\item[(a)] Koolen and Park \cite{shilla} showed the non-existence
of \drg s with the intersection array (iii) in Theorem
\ref{non-exist-7} and so did Juri\v{s}i\'{c} and Koolen
\cite{jurisic-koolen} for the intersection arrays (iv)-(vii).
\item[(b)] Suppose that $\Gamma$ is a \drg~with an intersection
array (i), (ii) or (iii) in Theorem \ref{non-exist-7}. By
\cite[Proposition 4.2.17]{bcn}, $\Gamma_3$(the graph with the
vertices are $V(\g)$ and the edges are the $2$-subsets of vertices
at distance $3$ in $\g$) is a strongly regular graph with
parameters $(672,121,20,22)$, $(855,126,21,18)$ or
$(924,143,22,22)$, respectively. No strongly regular graphs with
these parameters are known.
\end{enumerate}
\end{remark}

\section{Geometric \drg s with smallest eigenvalue $-3$}\label{geosection}
In this section, we prove Theorem \ref{gdrg} in which we determine the geometric distance-regular graphs with smallest eigenvalue $-3$.\\

Let $\Gamma$ be a distance-regular graph with diameter $D=D(\g)$. For any non-empty subset $X$ of $V(\Gamma)$ and for each $i=0,1,\ldots,D$, we put
\[
X_i:=\{x\in V(\Gamma) \mid d(x,X)=i\},
\]
where $d(x,X)=\min\{d(x,y) \mid y\in X \}$. Suppose that $C\subseteq V(\g)$ is a Delsarte clique in $\Gamma$. For each $i=0,1,\ldots, D-1$
and for a vertex $x \in C_i$, define
\[\psi_i(x,C) := \left| \{ z \in C \mid d(x,z) = i\}\right|.\]
The number $\psi_i(x,C)~(i=0,1,\ldots, D-1)$
depends not on the pair $(x,C)$ but depends only on the distance $i =
d(x,C)$ (cf. \cite[Section 4]{DCG1} and \cite[Section 11.7]{godsil-93}). Hence denote
\[\psi_i:=\psi_{i}(x,C)~~(i=0,1,\ldots,D-1).\]
Now, let $\Gamma$ be geometric with respect to $\mathcal{C}$ a set of Delsarte cliques. For $x, y\in V(\Gamma)$ with
$d(x,y)=i~~(i=1,2,\ldots, D)$, define $\tau_i(x,y;\mathcal{C}) $ as
the number of cliques $C$ in $\mathcal{C}$ satisfying $x\in C$ and $d(y,C) = i-1$. By \cite[Lemma 4.1]{DCG1}, the number
$\tau_i(x,y;\mathcal{C})$ ($i=1,2,\ldots,D$) depends not on the pair $(x,y)$ and $\mathcal{C}$, but depends
only on the distance $i=d(x,y)$. Thus we may
put
\[\tau_i:=\tau_i(x,y;\mathcal{C})~(i=1,2,\ldots,D)~.\]
Note that for any geometric distance-regular graph $\g$,
\begin{equation}\label{tau-D}
\tau_D=-\theta_{\min}
\end{equation}
holds, where $D=D(\g)$ and $\theta_{\min}=\theta_{\min}(\g)$.\\
The next lemma is a direct consequence of \cite[Proposition 4.2 (i)]{DCG1}.

\begin{lemma} \label{geo-para}
Let $\Gamma$ be a geometric \drg. Then the following hold.\\
(i) $b_i=-(\theta_{\min}+\tau_i)\left(1-\frac{k}{\theta_{\min}}-\psi_i\right)$ $(1\leq i\leq D-1)$.\\
(ii) $c_i=\tau_i\psi_{i-1}$  $(1\leq i\leq D)$.
\end{lemma}

Note that by (\ref{tau-D}) and Lemma \ref{geo-para} (ii), any geometric \drg~with diameter $D$ satisfies
\begin{equation}\label{cD}
c_D=(-\theta_{\min})\psi_{D-1}\geq -\theta_{\min}.
\end{equation}

\begin{lemma}\label{psi-bd}
Let $\g$ be a geometric \drg. Then
\begin{equation}\label{basic-ineq-psi1}
\psi_1\leq \tau_2\leq -\theta_{\min}.
\end{equation}
In particular, $\psi_1^2\leq c_2\leq \theta_{\min}^2$ holds.
\end{lemma}
\pf Let $x$ be a vertex and let $C$ be a Delsarte clique satisfying
$x \not \in C$. If there are two neighbors $y$ and $z$ of $x$ in
$C$, then two edges $(x,y)$ and $(x,z)$ lie in different Delsarte
cliques as $\g$ is geometric. This shows $\psi_1\leq \tau_2$. Note
that the number of Delsarte cliques containing any fixed vertex is
$-\theta_{\min}$, so that $\tau_i\leq -\theta_{\min}$ for all
$i=1,\ldots, D$. Hence, we find $\psi_1\leq \tau_2\leq
-\theta_{\min}$. In particular, it follows by Lemma \ref{geo-para}
(ii) and (\ref{basic-ineq-psi1}) that $\psi_1^2\leq \tau_2
\psi_1=c_2\leq \theta_{\min}^2$ holds. \epf

\begin{theorem}\label{gdrg}
Let $\g$ be a geometric \drg~with smallest eigenvalue $-3$. Then $\g$ satisfies one of the following.
\begin{enumerate}
\item[(i)] $k=3$ and $\g$ is one of the following graphs: the Heawood graph, the Pappus graph, Tutte's $8$-cage, the Desargues graph, Tutte's $12$-cage, the Foster graph, $K_{3,3}$, $H(3,2)$.
\item[(ii)] A Steiner graph $S_3(\alpha-3)$, i.e., a geometric strongly regular graph with parameters $\left(\frac{(2\alpha-3)(\alpha-2)}{3},3\alpha-9,\alpha,9  \right)$, where $\alpha\geq 6$ and $\alpha \equiv 0,2~~(\mbox{mod}~3)$.
\item[(iii)] A Latin square graph $LS_3(\alpha)$, i.e., a geometric strongly regular graph with parameters $(\alpha ^2,3(\alpha-1),\alpha,6)$, where $\alpha \geq 4$.
\item[(iv)] The generalized $2D$-gon of order $(s,2)$, where $(D,s)=(2,2),(2,4),(3,8)$.
\item[(v)] One of the two generalized hexagons of order $(2,2)$ with $\iota(\g)=\{6,4,4;1,1,3\}$.
\item[(vi)] A generalized octagon of order $(4,2)$ with $\iota(\g)=\{12,8,8,8;1,1,1,3\}$.
\item[(vii)] The Johnson graph $J(\alpha,3)$, where $\alpha \geq 6$.
\item[(viii)] $D=3$ and $\iota(\g)=\{3\alpha +3,2\alpha +2, \alpha +2-\beta ;1,2,3\beta\}$, where $\alpha \geq \beta \geq 1$.
\item[(ix)] The halved Foster graph with $\iota(\g)=\{6,4,2,1;1,1,4,6\}$.
\item[(x)] $D=\mbox{\em \HH}+2\geq 4$ and
$$(c_i,a_i,b_i)=\left\{ \begin{array}{ll}
(1,\alpha,2\alpha+2)& \mbox{ for }1\leq i\leq \mbox{\em \HH}\\
(2,2\alpha+\beta-1,\alpha-\beta+2) & \mbox{ for }i=\mbox{\em \HH} +1\\
(3\beta,3\alpha-3\beta+3,0) & \mbox{ for }i=\mbox{\em \HH} +2
\end{array}
\right., \mbox{~where~} \alpha \geq  \beta \geq 2.$$
\item[(xi)] $D=\mbox{\em \HH}+2\geq 3$ and $$(c_i,a_i,b_i)=\left\{ \begin{array}{ll}
(1,\alpha,2\alpha+2)& \mbox{ for }1\leq i\leq \mbox{\em \HH} \\
(1,\alpha+2\beta-2,2\alpha-2\beta+4) & \mbox{ for }i=\mbox{\em \HH} +1\\
(3\beta,3\alpha-3\beta+3,0) & \mbox{ for }i=\mbox{\em \HH} +2
\end{array}
\right., \mbox{~where~} \alpha\geq \beta \geq 2.$$
\item[(xii)] A distance-$2$ graph of a distance-biregular graph with vertices of valency $3$ and
$$(c_i,a_i,b_i)=\left\{ \begin{array}{ll}
(1,\alpha,2\alpha+2)& \mbox{ for }1\leq i\leq \mbox{\em \HH}\\
(1,\alpha+2,2\alpha)& \mbox{ for }i=\mbox{\em \HH} +1\\
(4,2\alpha-1,\alpha)& \mbox{ for }\mbox{\em \HH} +2\leq i\leq D-2\\
(4,2\alpha+\beta-3,\alpha-\beta+2) & \mbox{ for }i=D-1\\
(3\beta,3\alpha-3\beta+3,0) & \mbox{ for }i=D
\end{array}
\right., \mbox{~where~} \alpha\geq \beta \mbox{~and~} \beta\in \{2,3\}.$$
\end{enumerate}
\end{theorem}

\pf Let $\g$ be geometric with respect to $\mathcal{C}$. As $\theta_{\min}=-3$, we have $k\equiv 0$ (mod $3$). If $k=3$ then $\g$ satisfies (i) by \cite{k=3} (cf.\cite[Theorem 7.5.1]{bcn}). In the rest of the proof, we assume $k\geq 6$ and let $D=D(\g)$. We divide the proof into two cases, ({\bf Case 1: $c_2\geq 2$}) and ({\bf Case 2: $c_2 = 1$}).\\

\noindent {\bf Case 1: $c_2 \geq 2$}\\
By (\ref{basic-ineq-psi1}) with $\theta_{\min}=-3$, we find $\psi_1\in \{1,2,3\}$.\\

First suppose $\psi_1=1$, so that $\g$ is locally a disjoint union of three cliques of size $a_1+1$ and $k=3(a_1+1)$. By \cite[Theorem 3.1]{yamazaki}, $\g$ satisfies either ($c_2=2$ and $2\leq D \leq 3$) or ($c_2=3$ and $D=2$). If $c_2=2$ and $D=2$ then $-3$ is not the smallest eigenvalue of the matrix $L_1(\g)$ in (\ref{mtx-L}), which contradicts to $\theta_{\min}=-3$. If $c_2=2$ and $D =3$ then $\tau_2=2$ and $\tau_3=3$ by Lemma \ref{geo-para} (ii) and (\ref{tau-D}), respectively, and thus $(c_1,a_1,b_1)=(1,a_1,2a_1+2)$, $(c_2,a_2,b_2)=(2,2a_1-1+\psi_2,a_1+2-\psi_2)$ and $(c_3,a_3,0)=(3\psi_2,3a_1+3-3\psi_2,0)$ all hold by Lemma \ref{geo-para}. Now, $\g$ satisfies (viii). If $c_2=3$ and $D =2$, then $\Gamma$ is the generalized quadrangle of order $(s,2)$, where $s=2,4$ (cf. \cite[Theorem 6.5.1]{bcn} and \cite[Theorem 1]{rnp(t=2)}).\\

Next suppose $\psi_1=2$, so that $\tau_2\in \{2,3\}$, $b_1=\frac{2(k-3)}{3}$ and $c_2=2\tau_2$ all follow by (\ref{basic-ineq-psi1}) and Lemma \ref{geo-para}.
If $D \geq  3$ then $\g$ is the Johnson graph $J(\alpha,3)~(\alpha \geq 6)$ of diameter $3$ by \cite[Theorem 7.1]{-m} and \cite[Remark 2 (ii)]{DCG2}.
Now, we consider $D =2$. Then, $\tau_2=3$ by (\ref{tau-D}), and $\Gamma$ is a strongly regular graph with parameters $(a_1^2,3(a_1-1),a_1,6)$, where $a_1\geq 4$
as $k\geq 6$ and $\g$ is geometric. Hence, (iii) follows as $\g$ is the line graph of a $2-(3\alpha,3,1)$-transversal design,
where $\mathcal{C}$ and $V(\g)$ are the set of points and lines respectively (See Remark \ref{geo-rmk} (b)).\\

Finally, we consider $\psi_1=3$. Then $c_2=\tau_2 \psi_1=9$ holds by Lemma \ref{psi-bd}. From Lemma \ref{geo-para} (i) with $\theta_{\min}+\tau_2=0$, $D=2$ follows,
and thus $(c_1,a_1,b_1)=(1,a_1,2a_1-10)$ and $(c_2,a_2,b_2)=(9,3a_1-18,0)$. Since $\g$ is geometric,
$\g$ is a Steiner graph $S_3(\alpha-3)$ and $\g$ satisfies (ii), where the restriction on $a_1$ is obtained from $k\geq 6$ and
the fact that $|V(\g)|$ is a positive integer (See \cite[p.396]{neumaier-m} and Remark \ref{geo-rmk}). This completes the proof of {\bf Case 1}.\\

\noindent {\bf Case 2: $c_2 = 1$}\\
From the conditions $c_2=\tau_2 \psi_1=1$ and $\theta_{\min}=-3$, $\g$ is locally a disjoint union of three cliques of size $a_1+1$. If $a_1\leq 1$ then $k\in \{3,6\}$ follows from $|C|\in \{2,3\}$ for any Delsarte clique $C$ in $\g$. By \cite{a=1k=6}, $\g$ satisfies (v) or (ix). \\
From now on, we assume $a_1\geq 2$. First suppose $c_{\HH+1}\geq 2$, where $\HH=\HH(\g)$ is the head of $\g$ in (\ref{head}). Then by (\ref{cD}) and \cite[Theorem 3.1]{yamazaki}, $\g$ satisfies either ($c_{\HH+1}=3$ and $D=\HH+1$) or ($c_{\HH+1}=2$ and $D=\HH+2$). For the case $c_{\HH +1}=3$, $\g$ is a generalized $2D$-gon of order $(s,2)$, where $(D,s)=(3,8),(4,4)$ (cf. \cite[Section 6.5]{bcn} and \cite[Theorem 1]{rnp(t=2)}). If $c_{\HH+1}=2$, then we find $\psi_{\HH}=1$ and $\tau_{\HH+1}=2$ by $c_{\HH}=\psi_{\HH-1}\tau_{\HH}=1$ and \[a_1=a_{\HH}=\tau_{\HH}(a_1+1-\psi_{\HH-1})+(3-\tau_{\HH})(\psi_{\HH}-1),\]
from which (x) holds by (\ref{tau-D}), Lemma \ref{geo-para} and \cite[Proposition 2]{rnp(t=2)}. Next suppose $c_{\HH+1}=1$. By (\ref{cD}) and \cite[Theorem 4.1]{yamazaki}, $\g$ satisfies either $D=\HH+2$ or (xii). For the case $D=\HH+2$ with $c_{\HH+1}=1$, (xi) follows by (\ref{tau-D}) and Lemma \ref{geo-para}. This completes the proof of Theorem \ref{gdrg}. \epf

We remark on the \drg s in Theorem \ref{gdrg}.

\begin{remark}\label{geo-rmk}
\begin{enumerate}
\item[(a)] The line graph of a Steiner triple system on $2\alpha-3$ points for any integer $\alpha \geq 6$ satisfying $\alpha \equiv 0,2$ (mod~$3$), which is called a {\em Steiner graph} $S_3(\alpha-3)$, is a \srg~given in (ii). With the fact that a Steiner triple system on $v$ points exists for each integer $v$ satisfying $v\equiv 1$ or $3$ (mod $6$), Wilson showed in \cite{wilson74} and \cite{wilson75} that there are super-exponentially many Steiner triple systems for an admissible number of points, hence so are \srg s in (ii) (cf. \cite[p.~209]{cameronsrg}, \cite[Lemma 4.1]{neumaier-m}).
\item[(b)] The line graph of a $2-(mn,m,1)$-transversal design ($n\geq m+1$) is called a Latin square graph $LS_m(n)$ (See \cite[p.396]{neumaier-m}). In particular, a Latin square graph $LS_3(\alpha)$ is a geometric strongly regular graph in (iii). Since there are more than exponentially many Latin squares of order $\alpha$, so are such strongly regular graphs in (iii) (cf. \cite[p.~210]{cameronsrg}, \cite[Lemma 4.2]{neumaier-m}).
\item[(c)] In the list of \cite[Chapter 14]{bcn}, only the Hamming graph $H(3,\alpha+2)$, the Doob graph of diameter $3$ and the intersection array $\{45,30,7;1,2,27\}$ satisfy (viii). No \drg~with the last array, $\{45,30,7;1,2,27\}$, is
known. We can also check that if $\Gamma$ satisfies (viii) then the
eigenvalues of $\Gamma$ are integers.
\end{enumerate}
\end{remark}

\noindent {\em Proof of Theorem \ref{main-cor}: } It is straightforward from Theorem \ref{main-thm} and Theorem \ref{gdrg}. \epf

\begin{center}
{\bf Acknowledgements}
\end{center}
The author was supported by the Korea Research Foundation Grant
funded by the Korean Government(MOEHRD, Basic Research Promotion
Fund) KRF-2008-359-C00002. The author would like to thank Jack
Koolen for his valuable comments, and Jongyook Park for his careful
reading.

\end{document}